\documentclass[a4paper,11pt,reqno]{article}
\usepackage[T1]{fontenc}
\usepackage{courier}
\usepackage[scaled=0.92]{helvet}

\usepackage{amscd,amsmath,amssymb,amsthm,amsfonts,epsfig,graphics}

\usepackage[english]{babel}
\usepackage{amsthm,amsfonts,amssymb,amsmath}
\usepackage{hyperref}
\hypersetup{
colorlinks=true,
citecolor=red,
linkbordercolor={1 1 1},
citebordercolor={1 1 1},
}

\usepackage{mathrsfs}

\usepackage{draftwatermark}
\SetWatermarkLightness{1} %1=white, 0=black
\SetWatermarkAngle{45}
\SetWatermarkText{DRAFT}
\usepackage{xcolor} 

\theoremstyle{plain}
\newtheorem{satz}{Theorem}[section]
\newtheorem{prop}[satz]{Proposition}

\newtheorem{lem}[satz]{Lemma}

\theoremstyle{definition}

\newtheorem{rem}[satz]{Remark}

\newtheorem{hyp}[satz]{Hypothesis}

\newcommand{\mx}{\mbox}
\newcommand{\rw}{\rightarrow}
\newcommand{\de}{\displaystyle}
\newcommand{\ml}{\mathcal}
\newcommand{\tf}{\textbf}

\newcommand{\pl}{\partial}

\newcommand{\x}{\times}

\newcommand{\beq}[1]{\begin{equation} \label{#1}}
\newcommand{\eeq}{\end{equation}}
\newcommand{\beqar}{\[ \begin{array}{rcl}}
\newcommand{\eeqar}{\end{array} \]}

\providecommand{\ep}{\varepsilon}
\providecommand{\ph}{\varphi}

\providecommand{\RR}{\mathbb{R}}
\providecommand{\CC}{\mathbb{C}}

\providecommand{\ZZ}{\mathbb{Z}}
\providecommand{\NN}{\mathbb{N}}
\providecommand{\TT}{\mathbb{T}}

\newcommand{\anorm}[2]{\left \lVert#1 \right\rVert_{(1-d_{#2})(\rho,\sigma)}}
\newcommand{\cnorm}[2]{\left \lVert#1 \right\rVert_{(1-#2 d_j)(\rho_j,\sigma_j)}}
\newcommand{\snorm}[2]{\left| #1\right|_{#2}}
\newcommand{\norm}[2]{\left \lVert#1 \right\rVert_{#2}}

\newcommand{\lie}[1]{\mathcal{L}_{\chi^{(#1)}}}

\newcommand{\ah}{1/2}
\usepackage[utf8]{inputenc}

\hypersetup{pdfauthor={Fortunati, Wiggins},%
            pdftitle={Negligibility of small divisor effects...},%
            pdfsubject={Hamiltonian Mechanics},%
            pdfkeywords={Non-autonomous Hamiltonian systems, Aperiodic time dependence.}%
}

\DeclareMathOperator{\id}{Id}

\makeatletter
\renewcommand*{\@fnsymbol}[1]{\ensuremath{\ifcase#1\or *\or (a)\or (b)\or \else \fi}}
\makeatother
\title{\LARGE{\textbf{Negligibility of small divisor effects in the normal form theory for nearly-integrable Hamiltonians with decaying non-autonomous perturbations\footnote{This research was supported by ONR Grant No.~N00014-01-1-0769 and MINECO: ICMAT Severo Ochoa project SEV-2011-0087.}}}}
\date{}
\author{%
Alessandro Fortunati\thanks{E-mail: alessandro.fortunati@bristol.ac.uk} \and 
Stephen Wiggins\thanks{E-mail: s.wiggins@bristol.ac.uk}
\bigskip \\
School of Mathematics, University of Bristol, Bristol BS8 1TW, United Kingdom
}

\begin{document}

\maketitle

\begin{abstract}
The paper deals with the problem of the existence of a normal form for a  nearly-integrable real-analytic Hamiltonian with aperiodically time-dependent perturbation decaying (slowly) in time. In particular, in the case of an isochronous integrable part, the system can be cast in an exact normal form, regardless of the properties of the frequency vector.  The general case is treated by a suitable adaptation of the finite order normalization techniques usually used for Nekhoroshev arguments. The key point is that the so called ``geometric part'' is not necessary in this case. As a consequence, no hypotheses on the integrable part are required, apart from analyticity. \\
The work, based on two different perturbative approaches developed by A.Giorgilli et al., is a generalisation of the techniques used by the same authors to treat more specific aperiodically time-dependent problems.  
\smallskip\\
{\it Keywords:} Non-autonomous Hamiltonian systems, Aperiodic time dependence.
\smallskip\\
{\it 2010 MSC:} Primary: 70H08. Secondary: 37J40, 37J25.
 
\end{abstract}

\section{Introduction}
The problem of casting an analytic nearly-integrable Hamiltonian system into normal form is deeply related to Poincar\'{e}'s challenging \emph{probl\`{e}me g\'{e}n\'{e}ral de la dynamique} \cite{poi}. Nowadays, normal forms are still one of the main technical tools used to deal with the issue raised by Poincar\'{e} in this context.\\ 
The particular case in which the unperturbed part is supposed to be linear in the actions (isochronous case), already investigated by Birkhoff (and for this reason also known as the \emph{Birkhoff problem}) \cite{bir}, has a peculiar interest. The first rigorous statement concerning its stability can be found in  \cite{gal84}. The possibility to cast the considered Hamiltonian in normal form, up to some finite order\footnote{It is easy to see that any attempt to consider the limit $r \rw \infty$ would imply the degeneration into a trivial problem, (i.e. in which the allowed perturbation size reduces to zero, see also \cite[formula (46), Pag. 105]{gg85}).} $r$ and to obtain, as a consequence, a stability time estimate ``\`{a} la Nekhoroshev'', is directly related to a particularly simple small-divisors analysis: the non-resonant (Diophantine) hypothesis on the frequency vector $\omega$ of the unperturbed system is sufficient in order to ensure the resolvability of the (standard) homological equation arising in the normalization algorithm. An extensive bibliography on this problem goes beyond the purposes of this paper, we only mention the recent generalisations for the Planetary problem of \cite{pinz13} and of \cite{bamb} for infinite dimensional systems.\\
It is well known that the extension to the non-isochronous case requires a careful analysis (geometric part, see \cite{nek1}, \cite{nek2} and \cite{bg1986}) on the regions of the phase space in which the actions $I$ are such that $\omega=\omega(I)$ is non-resonant (non-resonant domains). \\  
The problem of dealing with time-dependent perturbations without any hypothesis on the time dependence (e.g. periodic or quasi-periodic) has peculiar technical difficulties. After the pioneering works of \cite{pust} and \cite{giozen}, the interest for this class of problems has been recently renewed in \cite{boun13},  \cite{fw14a} and subsequent papers. Basically, the novelty consists in the treatment of the time-dependent homological equation. A first approach consists in keeping the terms involving the time derivative of the generating function (also called extra-terms) in the normal form and then providing a bound for them. This approach, originally suggested in \cite{giozen} then used in \cite{fw14a}, yields a normal form result for the case a of slow time dependence. This hypothesis provides a smallness condition for the mentioned extra-terms. Alternatively, those terms can be removed by including them into the homological equation, which turns out to be, in this way, a linear ODE in time. This has been profitably used in \cite{fw14b}, \cite{fw15a} and in \cite{fw15} but requires (except for a particular case described in \cite{fw15}) an important assumption. More precisely, it is necessary to suppose that the perturbation, as a function of $t$, belongs to the class of summable functions over the real semi-axis\footnote{We stress that this hypothesis is usually not satisfied in the case of periodic or quasi-periodic time dependence.}. As in (\ref{eq:slowdecay}), those functions exhibiting a (slow) exponential decay will be used as a paradigmatic case. It will be shown that the consequences of this assumption in the isochronous case are remarkable:  the normalization algorithm can be iterated an infinite number of times by means of a \emph{superconvergent method} borrowed from KAM type arguments, see e.g. \cite{chierchia09}. The procedure leads to the so-called \emph{strong normal form} i.e. in which the normalized Hamiltonian has the same form of the integrable part of the initial problem. Furthermore, no restrictions are imposed  on $\omega$, hence flows with arbitrary frequencies persist in the transformed system. \\
As it would be likely to expect, this phenomenon has an important consequence also in the non-isochronous case. The possibility to disregard the problems related to the small divisors implies that the well known \emph{geography of the resonances} analysis, a key step of the Nekhoroshev theorem, is not necessary in this case and the results that can be stated are purely ``analytic''. In such a way, the classical assumptions on the unperturbed part of the Hamiltonian (such as steepness, convexity etc.), are no longer required. As a common feature with the isochronous case, the obtained normal form does not exhibit \emph{resonant terms}, as these  have been annihilated in the normalization by using the time-dependent homological equation. This implies that, in this case, the \emph{plane of fast drift} (see e.g. \cite{gio03}) degenerates to a point.\\
The paper uses in a concise but self-contained form, the tools developed in the above mentioned papers of the same authors, especially of \cite{fw15} in which the concept of ``family'' of canonical transformations parametrised by $t$ is introduced. The proofs are entirely constructed by using the language and the tools of the Lie series and Lie transform methods developed by Giorgilli et al., see e.g. \cite{gio03}. 
\section{Setting and main results} 
Consider the following nearly integrable Hamiltonian
\beq{eq:ham}
H(I,\ph,\eta,t)=h(I) + \eta + \hat{\ep} f(I,\ph,t) \mx{,}
\eeq
with $(I,\ph,\eta,t) \in G \times \TT^n \times \RR \times \RR^+$, where $ G \subset \RR^n$ and $\hat{\ep}>0$ is a small parameter, which is the ``autonomous equivalent'' in the extended phase space of Hamiltonian $\ml{H}(I,\ph,t)=h(I) + \hat{\ep} f(I,\ph,t)$.\\
We define, for all $t \in \RR^+:=[0,+\infty)$, the following complexified domain $\ml{D}_{\rho,\sigma}:=\ml{G}_{\rho} \times \TT_{\sigma}^n \x \ml{S}_{\rho}$, where $\ml{G}_{\rho}:=\bigcup_{I \in G} \Delta_{\rho}(I)$ and 
$$
\Delta_{\rho}(I):=\{\hat{I} \in \CC^n:|\hat{I}-I| \leq \rho\}, \quad \TT_{\sigma}^n := \{\ph \in \CC^n: |\Im \ph| \leq  \sigma\} \mx{,}
\quad
\ml{S}_{\rho}:=\{\eta \in \CC: |\Im \eta|  \leq \rho\} \mx{,}
$$
with $\rho,\sigma \in (0,1)$. For all $g:\mathcal{G}_{\rho} \times \TT_{\sigma}^n \times \RR^+ \rw \CC$, write $g=\sum_{k \in \ZZ^n} g_k(I,t) e^{i k \cdot \ph }$, then define the \emph{Fourier norm} (parametrized by $t$)
\beq{eq:fouriernorm}
\norm{g}{\rho,\sigma}:=\sum_{k \in \ZZ^n} \snorm{g_k(I,t)}{\rho} 
e^{|k|\sigma} \mx{,}
\eeq
with $|\cdot|_{\rho}$ is the usual supremum norm over $\ml{G}_{\rho}$ and 
$|k|:=\sum_{l=1}^n |k_l|$. For all $w : \mathcal{G}_{\rho} \times \TT_{\sigma}^n \times \RR^+ \rw \CC^n$ we shall set $\norm{w}{\rho,\sigma}:=\sum_{l=1}^n\norm{w_l}{\rho,\sigma}$  The standard framework (see eg. \cite{bggs84}) is the space $\mathfrak{C}_{\rho,\sigma}$, of continuous functions on $\ml{G}_{\rho} \times \TT_{\sigma}^n$, holomorphic in its interior for some $\rho,\sigma$ and real on $G \times \TT^n$ for all\footnote{In particular, if $g \in \mathfrak{C}_{\rho,\sigma}$ then $|g_k|_{\rho} \leq \norm{g}{\rho,\sigma} \exp(-|k|\sigma)$ for all $t \in \RR^+$.} $t \in \RR^+$. We shall suppose $h(I) \in \mathfrak{C}_{\rho,\cdot}$ and $f \in \mathfrak{C}_{\rho,\sigma}$ while it is sufficient to assume that, for all $I \in \ml{G}_{\rho}$, $f_k(I,\cdot) \in \ml{C}^1(\RR^+)$. \\
Similarly to \cite{fw15}, we introduce the following 
\begin{hyp}[Time decay] \label{hyp}
There exists $M_f>0$ and $a \in (0,1)$ 
\beq{eq:slowdecay}
\norm{f(I,\ph,t)}{\rho,\sigma} \leq M_f e^{-a t} \mx{.}
\eeq
\end{hyp}
Set $\ep:=\hat{\ep} M_f$. We firstly state the following
\begin{satz}[Strong aperiodic Birkhoff]\label{thm}
Consider Hamiltonian (\ref{eq:ham}) with $h(I):=\omega \cdot I$, under the Hypothesis \ref{hyp} and the described regularity assumptions. Then, for all $a \in (0,1)$ there exists $\ep_a>0$ such that the following statement holds true. For all $\ep \in (0,\ep_a]$, it is possible to find $0<\rho_*<\rho_0<\rho$ and $0<\sigma_*<\sigma_0<\sigma$ and an analytic, canonical, $\ep-$close and asymptotic to the identity change of variables  $(I,\ph,\eta) = \ml{B}(I^{(\infty)},\ph^{(\infty)},\eta^{(\infty)})$, $\ml{B}:\ml{D}_{\rho_*,\sigma_*} \rw \ml{D}_{\rho_0,\sigma_0}$ for all $t \in \RR^+$,  
casting Hamiltonian (\ref{eq:ham}) into the \emph{strong Birkhoff normal form} 
\beq{eq:birnormal}
H^{(\infty)}(I^{(\infty)},\ph^{(\infty)},\eta^{(\infty)})=\omega \cdot I^{(\infty)} + \eta^{(\infty)} \mx{.}
\eeq
\end{satz}
Hence, in the new variables, the flow with frequency $\omega$ persists \textbf{for all} $\omega$, regardless of the numerical features of this vector, i.e. more specifically, no matter if it is resonant or not. The absence of a non-resonance hypothesis on $\omega$ implies also that (\ref{eq:birnormal}) holds also if $\omega$ has an arbitrary number of zero components, implying the persistence of any \emph{lower dimensional torus}.
\\With a straightforward adaptation of the notational setting, the result in the general case states as follows:
\begin{satz}\label{thmtwo} There exist $\ep_a^*>0$ and $r \in \NN \setminus\{0\}$ such that, for all $\ep \in (0,\ep_a^*]$ it is possible to find an analytic, canonical, $\ep-$close and asymptotic to the identity change of variables  $(I,\ph,\eta) = \ml{N}_r(I^{(r)},\ph^{(r)},\eta^{(r)})$, $\ml{N}_r:\ml{D}_{\tilde{\rho_*},\tilde{\sigma_*}} \rw \ml{D}_{\tilde{\rho_0},\tilde{\sigma_0}}$ for all $t \in \RR^+$,  
casting Hamiltonian (\ref{eq:ham}) under the Hypothesis \ref{hyp}, into the \emph{normal form of order} $r$
\beq{eq:annormalform}
H^{(r)}(I^{(r)},\ph^{(r)},\eta^{(r)},t)=h(I^{(r)}) + \eta^{(r)} + \ml{R}^{(r+1)}(I^{(r)},\ph^{(r)},t) \mx{,}
\eeq
where $\ml{R}^{(r+1)}$ is ``exponentially small'' with respect to $r$ and vanishes for\footnote{See bound (\ref{eq:estrem}).} $t \rw +\infty$.
Moreover, for all $I(0) \in G$ one has in (\ref{eq:ham}): $|I(t)-I(0)| \leq \sqrt{\ep} \tilde{\rho}_0/8$ for all $t \in \RR^+$.
\end{satz}
Similarly to \cite{fw15} (and the mentioned previous papers), no lower bounds are imposed on $a$ so that the decay can be arbitrary slow. The (natural) consequence is that either $\ep_a$ or $\ep_a^*$ decrease with $a$, see (\ref{eq:epz}) and (\ref{eq:r}). 

\part{Proof of Theorem \ref{thm}}

\section{The normalization algorithm}
Given a function $G:=G (I,\ph,t)$, define the \emph{Lie series operator} $\exp({\ml{L}}_{G}):=\id+\sum_{s \geq 1} (1/s!) \ml{L}_{G}^s$, where $\ml{L}_{G} F :=\{F,G\} \equiv F_{\ph} \cdot G_I - G_{\ph} \cdot F_I- F_{\eta} G_t $. The aim is to construct a \emph{generating sequence} $\{\chi^{(j)}\}_{j \in \NN}$, such that the formal limit 
\beq{eq:limit}
\ml{B}:=\lim_{j \rw \infty} \ml{B}^{(j)} \circ \ml{B}^{(j-1)} \circ \ldots \circ \ml{B}^{(0)} \mx{,}
\eeq
where $\ml{B}^{(j)}:=\exp(\lie{j})$ is such that $\ml{B} \circ H$ is of the form (\ref{eq:birnormal}). The following statement shows that this is possible, at least at a formal level
\begin{prop}\label{prop:formal} Suppose that for some $j \in \NN$ Hamiltonian (\ref{eq:ham}) is of the form 
\beq{eq:hamric}
H^{(j)}=\omega \cdot I + \eta + F^{(j)} (I,\ph,t) \mx{.}
\eeq
Then $H^{(j+1)} :=\ml{B}^{(j)} \circ H^{(j)} $ is still of the form (\ref{eq:hamric}) with 
\beq{eq:fjp}
F^{(j+1)}=\sum_{s \geq 1} \frac{s}{(s+1)!}\lie{j}^s F^{(j)} \mx{,}
\eeq
provided that $\chi^{(j)}$ solves the \emph{homological equation}
\beq{eq:hom}
\chi_t^{(j)} + \omega \cdot  \chi_{\ph}^{(j)}=F^{(j)} \mx{.}
\eeq
\end{prop}
Since Hamiltonian (\ref{eq:ham}) is of the form (\ref{eq:hamric}), one can set $H^{(0)}:=H$ with $F^{(0)}:=\hat{\ep} f$. Thus, by induction, the form (\ref{eq:hamric}) holds for all $j \in \NN$. Clearly, this does not guarantee that the objects involved in the algorithm are meaningful for all $j$, as it is well known their sizes can grow unboundedly as $j$ increases,  as a consequence of \emph{small divisors} phenomena. The aim of Section \ref{sec:convergence} (and in particular of Lemma \ref{lem:main}) is to show that this is not the case: the key ingredient is the time decay of $f$.
\proof
We get $\exp(\lie{j})H^{(j)}=I \cdot \omega + \eta + F^{(j)}(I,\ph,t)+\lie{j}( \omega  \cdot I + \eta)+\sum_{s \geq 1} (1/s!) \lie{j}^s F^{(j)}+ \sum_{s \geq 2} (1/s!) \lie{j}^s (\omega  \cdot I + \eta)$. The sum between the third and fourth terms of the r.h.s. of the latter equation vanishes due to (\ref{eq:hom}). As for the last two terms, by setting $F^{(j+1)}$ as the sum of them, one gets $F^{(j+1)}=\sum_{s \geq 1} (1/s!) \lie{j}[F^{(j)}+(s+1)^{-1}\lie{j}(\omega  \cdot I +\eta)]$, which immediately yields (\ref{eq:fjp}) by using (\ref{eq:hom}).   
\endproof
The (formal) expansions $\chi^{(j)}=\sum_{k \in \ZZ^n} c_k^{(j)}(I,t)e^{i k \cdot \ph}$ and $F^{(j)}=\sum_{k \in \ZZ^n} f_k^{(j)}(I,t)e^{i k \cdot \ph}$ yield (\ref{eq:hom}) in terms of Fourier components
\beq{eq:homcomp}
\pl_t c_k^{(j)}(I,t)+i \lambda(k) c_k^{(j)} (I,t)= f_k^{(j)}(I,t) \mx{,} 
\eeq
with $\lambda(k):=\omega \cdot k$. The solution of (\ref{eq:homcomp}) is  
\beq{eq:solhom}
c_k^{(j)}(I,t)=e^{-i \lambda(k) t}\left[ c_k^{(j)}(I,0)+\int_0^t e^{i \lambda(k) s} f_k^{(j)}(I,s) ds\right] \mx{,}
\eeq 
where $ c_k^{(j)}(I,0)$ will be chosen later.
\section{Convergence}\label{sec:convergence}
The classical argument requires the construction of a sequence of nested domains $\ml{D}_{\rho_{j+1},\sigma_{j+1}} \subset \ml{D}_{\rho_j,\sigma_j} \ni (I^{(j)},\ph^{(j)},\eta^{(j)})$, such that $\ml{B}_j: \ml{D}_{j+1} \rw \ml{D}_{j}$. The resulting progressive restriction is essential in order to use standard Cauchy tools, see Prop. \ref{prop:poisson}. The estimates found in Lemma \ref{lem:homeq}, concerning the solution of equation (\ref{eq:hom}), will be used to prove Lemma \ref{lem:main}, providing in this way the bound on $F^{(j)}$ defined in Prop. \ref{prop:formal}. This is achieved for a suitable sequence of domains prepared in Lemma \ref{lem:preparation} via $\{\rho_j\}$ and $\{\sigma_j\}$. This allows us to conclude that the perturbation term is actually removed in the limit (\ref{eq:limit}).\\
The final step consists of showing that $\ml{B}$ defines an analytic map $\ml{B} :  \ml{D}_{\rho_*,\sigma_*} \ni (I^{(\infty)},\ph^{(\infty)},\eta^{(\infty)}) \rw  \ml{D}_{\rho_0,\sigma_0} \ni (I^{(0)},\ph^{(0)},\eta^{(0)}) \equiv (I,\ph,\eta)$, where $\rho_* \leq \rho_j$ and $\sigma_* \leq \sigma_j$ for all $j \in \NN$. This property is shown in Lemma \ref{lem:transf}. As $\ml{D}_{\rho_*,\sigma_*}$ will be the domain of analyticity of the transformed Hamiltonian via $\ml{B}$, it will be essential to require  that $\rho_*,\sigma_*>0$. 
\subsection{Some preliminary results}
\begin{prop}\label{prop:poisson}
Let $F,G: \mathcal{G}_{\rho} \times \TT_{\sigma}^n \times \RR^+ \rw \CC$ such that $\norm{F}{(1-d')(\rho,\sigma)}$ and $\norm{G}{(1-d'')(\rho,\sigma)}$ are bounded for some $d',d'' \in [0,1)$. Then, defining $\delta:=|d'-d''|$ and $\hat{d}:=\max\{d',d''\}$, for all $\tilde{d} \in (0,1-\hat{d})$ one has for all $s \in \NN \setminus \{0\}$
\beq{eq:lie}
\norm{\ml{L}_G^s F}{(1-\tilde{d}-\hat{d})(\rho,\sigma)} \leq  
\frac{s!}{e^2} \left( \frac{2e}{\tilde{d}(\tilde{d}+\tilde{\delta}_s)\rho \sigma} \norm{G}{(1-d'')(\rho,\sigma)} \right)^s
\norm{F}{(1-d')(\rho,\sigma)} \mx{,}
\eeq
where $\tilde{\delta}_s=\delta$ if $s=1$ and is zero otherwise.
\end{prop}
\proof
Straightforward from \cite[Lemmas 4.1, 4.2]{gio03}.
\endproof
\begin{lem}\label{lem:homeq}
Suppose that $F^{(j)}$ satisfies $\norm{F^{(j)}}{[\hat{\sigma},\hat{\rho}]} \leq M^{(j)}\exp(-at)$ for some $M^{(j)}>0$, $\hat{\rho} \leq \rho$ and $\hat{\sigma} \leq \sigma$. Define $C_{\omega}:=1+|\omega|$, then for all $\delta \in (0,1)$ the solution of (\ref{eq:hom}) satisfies
\beq{eq:boundshom}
\norm{\chi^{(j)}}{(1-\delta)(\hat{\rho},\hat{\sigma})} \leq 
\frac{M^{(j)}}{a} \left( \frac{e}{\delta \hat{\sigma}} \right)^{2n} e^{-at}, \qquad 
\norm{\chi_t^{(j)}}{(1-\delta)(\hat{\rho},\hat{\sigma})} \leq 
C_{\omega} \frac{M^{(j)}}{a} \left( \frac{e}{\delta \hat{\sigma}} \right)^{2n} e^{-at} \mx{.} 
\eeq  
\end{lem}
\proof First of all, by hypothesis $|f_k^{(j)}(I,t)| \leq M^{(j)} \exp(-|k|\hat{\sigma}-at)$, in particular, by choosing $c_k^{(j)}(I,0):=-\int_{\RR^+} \exp(i \lambda(k)s) f_k^{(j)}(I,s) ds$ we have that $|c_k^{(j)}(I,0)|<+\infty$ for all $I \in \ml{G}_{\rho}$. Substituting $c_k^{(j)}(I,0)$ in (\ref{eq:solhom}) one gets $|c_k^{(j)}(I,t)| \leq \int_t^{\infty}|f_k^{(j)}(I,s)|ds \leq (M^{(j)} /a)  \exp(-|k|\hat{\sigma} -at)$ which yields\footnote{Recall (\ref{eq:fouriernorm}), then use the inequality $\sum_{k \in \ZZ^n} \exp(-\delta|k|\hat{\sigma}) \leq (e \delta^{-1} \hat{\sigma}^{-1})^{2n}$. Its variant $\sum_{k \in \ZZ^n} (1+|\omega||k|)\exp(-\delta|k|\hat{\sigma}) \leq C_{\omega} (e \delta^{-1} \hat{\sigma}^{-1})^{2n}$ is used to obtain the second of (\ref{eq:boundshom}).} the first of (\ref{eq:boundshom}). As for the second of (\ref{eq:boundshom}), it is sufficient to use (\ref{eq:homcomp}), which implies, $|\pl_t c_k^{(j)}(I,t)| \leq (M^{(j)}/a)(1+|\omega||k|) \exp(-|k|\hat{\sigma}-a t)$ then proceed similarly.
\endproof
\begin{rem}
It is immediate to notice that a hypothesis of non-resonance on $\omega$ does not substantially improve the bounds (\ref{eq:boundshom}). A more careful computation yields 
$$
|c_k^{(j)}(I,t)| \leq M^{(j)} (a^2+(\omega \cdot k)^2)^{-\frac{1}{2}}e^{- |k| \sigma_j-a t} \mx{,}
$$ 
Hence the estimate cannot be refined due to the presence of $|c_0^{(j)}(I,t)|$, no matter what the minimum value of $(\omega \cdot k)$ is.
\end{rem}
\subsection{A suitable sequence of domains}
\begin{lem}\label{lem:preparation}
Let $\{d_j\}_{j \in \NN}$ be a (real valued) sequence such that $0 \leq d_j \leq 1/6$. Consider, for all $j \in \NN$, the following sequences 
\beq{eq:iterative}
\epsilon_{j+1}:=K a^{-1} d_j^{-\tau} \epsilon_j^2, \qquad (\rho_{j+1},\sigma_{j+1}):=(1-3d_j)(\rho_j,\sigma_j) \mx{,}
\eeq
with $K>0$ and $\tau:=2 n + 3$. Then, for all $0<\rho_0 \leq \rho$, $0<\sigma_0 \leq \sigma$ and $\epsilon_0 \leq \ep_a$ where 
\beq{eq:epz}
\ep_a \leq a K^{-1} (2 \pi)^{-2 \tau} \mx{,}
\eeq
it is possible to construct $\{d_j\}_{j \in \NN}$ such that $(\rho_*,\sigma_*)=(1/2)(\rho_0,\sigma_0)$, in particular they are strictly positive. Furthermore $\lim_{j \rw \infty} \epsilon_j=0$. 
\end{lem}
\proof Choose $\epsilon_j:=\epsilon_0 (j+1)^{-2 \tau}$ (so that $\lim_{j \rw \infty} \epsilon_j=0$ by construction). By the first of (\ref{eq:iterative}) one gets 
\beq{eq:dj}
d_j = (\epsilon_0 K a^{-1})^{\frac{1}{\tau}} (j+2)^2/(j+1)^4 \mx{,}
\eeq
hence, by (\ref{eq:epz}), $d_j \leq \pi^{-2} (j+1)^{-2}$. This implies $\sum_{j \geq 0} d_j \leq 1/6$ and then, trivially, $d_j \leq 1/6$ for all $j \in \NN$. Now we have\footnote{Use the inequality $ \ln (1-x) \geq -2 x \ln 2$, valid for all $x \in [0,1/2]$.} $ \ln \Pi_{j \geq 0} (1-3 d_j)=\sum_{j \geq 0} \ln (1-3 d_j) \geq -6 \ln 2 \sum_{j \geq 0} d_j = - \ln 2$, hence $\lim_{j \rw \infty} \rho_j=\rho_0 \Pi_{j \geq 0} (1-3 d_j) \geq \rho_0/2 =: \rho_*$. Analogously $\sigma_*:=\sigma_0/2$.
\endproof
\subsection{Bounds on the formal algorithm}
\begin{lem}\label{lem:main} There exists $K=K(\rho_0,\sigma_0)>0$ such that, if $\ep \leq  \ep_a$ where $\ep_a$ satisfies (\ref{eq:epz}), then
\beq{eq:decay} 
\norm{F^{(j)}}{(\rho_j,\sigma_j)}\leq \epsilon_j e^{-at} \mx{,}
\eeq
for all $j \in \NN$. Hence, the transformed Hamiltonian $\ml{B} \circ H$ is in the form (\ref{eq:birnormal}).
\end{lem}  
\proof By induction. Note that (\ref{eq:decay}) is true for $j=0$ setting $\epsilon_0:=\ep$. The condition on $\ep$ ensures the validity of Lemma \ref{lem:preparation}. Hence, supposing (\ref{eq:decay}), by Lemma \ref{lem:homeq} and Lemma \ref{lem:preparation}, we get
\beq{eq:boundchi}
\cnorm{\chi^{(j)}}{} \leq \epsilon_j (e/\sigma_*)^{2n} a^{-1} d_j^{-2n} e^{-a t} \mx{.}
\eeq
By (\ref{eq:fjp}) and Prop. \ref{prop:poisson} with $d'=d_j$, $d''=0$ and $\tilde{d}=d_j$ (the condition $d_j \leq 1-d_j$ holds as $d_j \leq 1/6$)
\beq{eq:fjpfirst}
\cnorm{F^{(j+1)}}{2} \leq \sum_{s \geq 1} \frac {1}{s!} \cnorm{\lie{j}^s F^{(j)}}{2} \leq 2^{-1} \Theta \norm{F^{(j)}}{(\rho_j,\sigma_j)} \mx{,}
\eeq
where\footnote{The reason for using $n C_{\omega}$ in the definition of $\Theta$ will be clear in the proof of Lemma \ref{lem:transf}.} 
\beq{eq:convergenceone}
\Theta:=2 \epsilon_j n C_{\omega} (e/\sigma_*)^{\tau} \rho_*^{-1}
a^{-1} d_j^{-2n-2} e^{- a t} 
\leq 1/2
\eeq
is a sufficient condition for the convergence of the operator $\exp(\lie{j})$, from which $\sum_{s \geq 1} \Theta^s \leq 2 \Theta$. Hence, by (\ref{eq:fjpfirst}), (\ref{eq:convergenceone}), then by (\ref{eq:boundchi}) one gets (use also $\sigma_*,\rho_*,d_j<1$)
\beq{eq:fjpsecond}
\cnorm{F^{(j+1)}}{2} \leq \epsilon_j^2 n C_{\omega}
 (e/\sigma_*)^{\tau}  \rho_*^{-1} a^{-1} d_j^{-\tau} e^{-at} \mx{.}
\eeq
The latter is valid \emph{a fortiori} in $\ml{D}_{(1-3d_j)(\rho_j,\sigma_j)}$.\\
In conclusion, by choosing $K:=n C_{\omega} (e/\sigma_*)^{\tau} \rho_*^{-1} =2^{\tau+1} n C_{\omega} (e/\sigma_0)^{\tau} \rho_0^{-1}$, from the first of (\ref{eq:iterative}), we have that (\ref{eq:decay}) is satisfied for $j \rw j+1$. Furthermore, by the first of (\ref{eq:iterative}), condition (\ref{eq:convergenceone}) yields $1 \geq 4 \epsilon_j K d_j a^{-1} d_j^{-\tau}e^{-at} =4 d_j (\epsilon_{j+1}/\epsilon_j) e^{-a t}$. The latter is trivially true for all $t \in \RR^+$ by the monotonicity of $\epsilon_j$ and as $d_j \leq 1/6$. Furthermore this implies
\beq{eq:nabla}
\Theta \leq 2 d_j e^{-at} \mx{.}
\eeq
Hence $\exp(\lie{j})$ is well defined for all $j \in \NN$. 
\endproof 
In this way the value of $\ep_a$ mentioned in the statement of Theorem \ref{thm} is determined once and for all. \\
\subsection{Estimates on the transformation of coordinates}
\begin{lem}\label{lem:transf}
The limit (\ref{eq:limit}) exists, it is $\ep-$close to the identity and satisfies
\beq{eq:boundtransf}
|I^{(\infty)}-I|,|\eta^{(\infty)}-\eta| \leq (\rho_0/6)e^{-at}\mx{,} \qquad |\ph^{(\infty)}-\ph| \leq (\sigma_0/6)e^{-at} \mx{,}
\eeq
in particular it defines an analytic map  $\ml{B}:\ml{D}_{\rho_*,\sigma_*} \rw \ml{D}_{\rho_0,\sigma_0}$ and $H^{(\infty)}$ is an analytic function on $\ml{D}_{\rho_*,\sigma_*} $ for all $t \in \RR^+$.
\end{lem}
\proof
Let us start with $I$. Note that $\norm{\lie{j} I^{(j+1)}}{(1-2d_j)(\rho_j,\sigma_j)} \leq n (e d_j \rho_j)^{-1}\cnorm{\chi^{(j)}}{}$ by a Cauchy estimate (see \cite[Lemma 4.1]{gio03}), so that the presence of $n$ in (\ref{eq:convergenceone}) is justified. Hence use Prop. \ref{prop:poisson} with $F \leftarrow \lie{j} I^{(j+1)}$, $s \leftarrow s-1$, obtaining
$\norm{\lie{j}^{s} \ph^{(j+1)}}{(1-3d_j)(\rho_j,\sigma_j)} \leq e^{-2}s!\Theta^{s} \rho_0$. This implies 
\[
|I^{(j+1)}-I^{(j)}|\leq e^{-2}\sum_{s \geq 1} (1/s!) \norm{\lie{j}^{s} I^{(j+1)}}{(1-3d_j)(\rho_j,\sigma_j)} \leq 2^{-1} \Theta \rho_0 \leq d_j \rho_0 e^{-at} \mx{,}
\]  
by (\ref{eq:nabla}). In particular $|I^{(j+1)}-I^{(j)}|$ is $\ep-$ close to the identity by (\ref{eq:dj}) for all $j \in \NN$, hence $|I^{(\infty)}-I| \leq \sum_{j \geq 0} |I^{(j+1)}-I^{(j)}|$ is. It is now sufficient to recall $\sum_{j \geq 0} d_j \leq 1/6$ in order to conclude.\\ 
The argument for $\ph$ is analogous while the variable $\eta$ requires a slight modification. In particular, as one needs to set $F \leftarrow \lie{j} \eta=-\chi_t^{(j)}$, the use of the second of (\ref{eq:boundshom}) requires the contribution of $C_{\omega}$ in (\ref{eq:convergenceone}). \\
In conclusion, the obtained composition of analytic maps is uniformly convergent in any compact subset of $\ml{D}_{\rho_*,\sigma_*}$. This implies that $\ml{B}$ is analytic on $\ml{D}_{\rho_*,\sigma_*}$ by the Weierstra{\ss} Theorem and hence the image of $H$ via $\ml{B}$ is an analytic function in the same domain.  
\endproof

\section{Further perturbation examples}
In this section we consider two alternative examples of perturbation. The main purpose is to show that the hypothesis of summability in time over the semi-axis is the only key requirement for the argument beyond the proof of  Theorem \ref{thm}.\\ 
In particular, we shall firstly consider a decay which is assumed to be quadratic in time, while in the second example a perturbation exhibiting a finite number of (differentiable) bumps is examined. The procedure is fully similar, with the exception of some bounds that will be explicitly given below.
\subsection{Quadratic decay}
Let us suppose that (\ref{eq:slowdecay}) is modified as 
$$
\norm{f(I,\ph,t)}{\rho,\sigma} \leq M_f (t+1)^{-2}
\mx{.}
$$
In the same framework, it is immediate to show that the analogous of Lemma \ref{lem:homeq} yields the following estimates
\[
\norm{\chi^{(j)}}{(1-\delta)(\hat{\rho},\hat{\sigma})} \leq M^{(j)} (e \delta^{-1} \hat{\sigma}^{-1})^{2n}  (t+1)^{-1}, \qquad \norm{\chi_t^{(j)}}{(1-\delta)(\hat{\rho},\hat{\sigma})} \leq M^{(j)} C_{\omega} (e \delta^{-1} \hat{\sigma}^{-1})^{2n} (t+1)^{-1} \mx{.} 
\]  
Clearly, in this case, the integration has led to a ``loss of a power'' in the decay. This is harmless as, by (\ref{eq:fjpfirst}), $\norm{F^{(j+1)}}{(1-2d_j)(\rho_j,\sigma_j)}=O(F^{j}) O(\chi^{(j)})+h.o.t.$ and then $F^{(j+1)} \sim (t+1)^{-3}  \leq (t+1)^{-2}$ so that the scheme can be iterated\footnote{A similar (and even stronger) phenomenon could have been noticed in the original setting. Namely, suppose by induction that $\norm{F^{(j)}}{(\rho_j,\sigma_j)} \leq \epsilon_j \exp(-a_j t)$. By Lemma \ref{lem:homeq} and  (\ref{eq:fjpfirst}), one finds that $\norm{F^{(j+1)}}{(\rho_{j+1},\sigma_{j+1})} \leq \epsilon_{j+1} \exp(- 2 a_j t)$ and so on. This leads to a remarkable rate of decay ($a_j=2^{j}a$) but not to a substantial improvement of the estimates and of the threshold (\ref{eq:epz}) of $\ep_a$, as these are uniform in $j$.}.\\
The rest of the proof is analogous provided that the term $e^{-at}$ is replaced with $1$ in the remaining estimates.
\subsection{Differentiable bumps} 
Let $L \in \NN \setminus \{0\}$ and $h>0$. Consider an increasing sequence $\{t_l\}_{l=1,\ldots,L} \in \RR^+$ such that $t_{l+1}-t_l>2h$, then the following function 
$$
\xi_l(t):= 
\left\{
\begin{array}{lcl}
(a_l/h^4)[(t-t_l+h)(t-t_l-h)]^2 & \quad &  t \in [t_l-h,t_l+h]\\
0 & & \mbox{otherwise} 
\end{array}
\right.
$$
where $a_l \in \RR$. Considering a function $\tilde{f}(I,\ph) \in \mathfrak{C}_{\rho,\sigma}$, we set as 
\[
f(I,\ph,t) := \tilde{f}(I,\ph) \sum_{l =1}^L \xi_l(t) \mx{.}
\]
In such case we find
\[
\norm{\chi^{(j)}}{(1-\delta)(\hat{\rho},\hat{\sigma})} \leq 2 A M^{(j)} h 
 (e \delta^{-1} \hat{\sigma}^{-1})^{2n}  , 
\qquad 
\norm{\chi_t^{(j)}}{(1-\delta)(\hat{\rho},\hat{\sigma})} \leq M^{(j)} C_{\omega} (e \delta^{-1} \hat{\sigma}^{-1})^{2n}\mx{,} 
\]
with $A:=\sum_{l=1}^L |a_l|$. The remaining part of the proof is straightforward with the obvious modifications.  In particular, as for the proof of Lemma \ref{lem:main}, one finds $K=2 n C_{\omega} (e/\sigma_*)^{\tau} h A \rho_*^{-1}$. 

\part{Proof of Theorem \ref{thmtwo}}
In order to simplify the notation, we shall use $(\rho_H,\sigma_H)$ in place of $(\rho,\sigma)$ and $(\rho,\sigma)$ in place of $(\tilde{\rho}_0,\tilde{\sigma}_0)$ from now on.
\section{Formal algorithm}
As in \cite{gio03}, we write Hamiltonian (\ref{eq:ham}) in the form
\[
H(I,\ph,\eta,t)=H_0(I,\eta)+H_1(I,\ph,t)+H_2(I,\ph,t)+\ldots
\]
where 
\[
H_0(I,\eta):=h(I)+\eta, \qquad H_s(I,\ph,t):=\sum_{k \in \Lambda_s} 
f_k(I,t)e^{i k \cdot \ph} \mx{,}
\] 
where $\Lambda_s:=\{k \in \ZZ^n: (s-1)N \leq |k| < sN\}$ and $N \in \NN \setminus \{0\}$ is meant to be determined.\\
Given a sequence of functions $\{\chi^{(s)}\}_{s\geq 1}:\mathfrak{C}_{\rho,\sigma} \rw \CC$, the Lie transform operator is defined as 
\beq{eq:lietransform}
T_{\chi}:=\de \sum_{s \geq 0} E_s,\qquad 
E_s:=
\left\{
\begin{array}{lcl}
\id & \quad& s=0\\
\de \frac{1}{s} \sum_{j=1}^s j \ml{L}_{\chi^{(j)}} E_{s-j} & \quad & s \geq 1
\end{array}
\right. \mx{.}
\eeq
Let $r \in \NN \setminus \{0\}$ to be determined. A \emph{finite} generating sequence of order $r$, denoted with $\chi^{[r]}$, is such that $\chi^{(s)} \equiv 0$ for all $s >r$. Our aim is to determine it in such a way the effect of $H_1,\ldots,H_r$ is removed, i.e. 
\beq{eq:normal}
H^{(r)}:=T_{\chi^{[r]}} H=H_0+\ml{R}^{(r+1)}(I,\ph,t) \mx{,}
\eeq
where the \emph{remainder} $\ml{R}^{(r+1)}$ contains $H_{> r}$ and a moltitude of terms produced during the normalization, which Fourier harmonics lie on $\Lambda_{>r}$. The smallness of the remainder is an immediate consequence of the decay property of the coefficients of an analytic function. The procedure is standard: condition (\ref{eq:normal}), with the use of (\ref{eq:lietransform}), yields a well known diagram which $s-$th \emph{level}\footnote{Namely, those terms of the diagram which Fourier harmonics belong to $\Lambda_s$.} is of the form 
\beq{eq:diagram}
\ml{E}_s:=E_s H_0+\sum_{l=1}^{s-1} E_{s-l} H_l + H_s=0 \mx{,}
\eeq 
if $s  = 2,\ldots,r$ and $E_1 H_0+H_1=0$ if $s=1$. As sum of all the ``non-normalised'' levels, the remainder easily reads as
\beq{eq:remaindernekho} 
\ml{R}^{(r+1)}=\sum_{s > r} \ml{E}_s \mx{.}
\eeq
By writing the first term of (\ref{eq:diagram}) in the form $E_s=\lie{s}+\sum_{j=1}^{s-1} (j/s) \lie{j}E_{s-j}$ and using the manipulation described in \cite[Chapter 5]{gio03}, one obtains a remarkable cancellation of the contribution of $H_0$. In this way, the generating sequence is determined as a solution of 
\beq{eq:homological}
\ml{L}_{H_0} \chi^{(s)} = \Psi_s,\qquad 
\Psi_s:=
\left\{
\begin{array}{lcl}
H_1 & \quad & s=1\\
\displaystyle H_s+\sum_{j=1}^{s-1} \frac{j}{s} E_{s-j}H_j  & & s \geq 2
\end{array}
\right. \mx{.}
\eeq
A formal expansion of $\chi^{(j)}$ and of $\Psi_s:=\sum_{k \in \ZZ^n} \psi_k^{(s)}(I,t)e^{i k \cdot \ph}$ yields for all $s=1,\ldots,r$ 
\beq{eq:homeqan}
\pl_t c_k^{(s)}(I,t)+i(\omega(I)\cdot k) c_k^{(s)}(I,t) = \psi_k^{(s)}(I,t), \qquad k \in \Lambda_s \mx{,}
\eeq
where, as usual, $\omega(I):=\pl_I h(I)$.
\begin{rem}
As a substantial difference with the isochronous case, the function 
$\omega(I)$ is a complex valued vector as $I \in \ml{G}_{\rho}$. In this way the exponent $\lambda(k)t$ appearing in formula (\ref{eq:solhom}) is no longer purely complex. More precisely, one finds a term of the form $\exp((\omega_C (I)\cdot k) t)$, having denoted $\omega(I)=\omega_R(I)+i \omega_C(I)$, $\omega_{R,C}(I) \in \RR^n$. The size of this term cannot be controlled without a cut-off on $k$. By restricting the analysis on the levels $\Lambda_s$ and using the fact that $|\omega_C(I)| \rw 0 $ as $ \rho \rw 0$, a loss ``of part of time decay'' at each step (see Lemma \ref{lem:homologicaleq}) will be the key ingredient to overcome this difficulty. The mentioned elements are clear obstructions to the limit $r \rw \infty$.    
\end{rem}

\section{Convergence}
\subsection{Set-up and some preliminary results}
The use of the analytic tools requires the usual construction of a sequence of nested domains. We shall choose, for all $s=1,\ldots,r$, the rule 
\beq{eq:domrestr}
d_s:=d (s-1)/r \mx{, }
\eeq
with $d \in (0,1/4]$. Clearly $d_s < d$ for all $s=1,\ldots,r$. Consider also the monotonically decreasing sequence of non-negative real numbers $\{a_s\}$ defined as follows 
\beq{eq:seqa}
a_{s+1}:=a_s(2r-s)/(2r), \qquad a_1:=a \mx{.}
\eeq
Given the analyticity domain of $H$ expressed by $(\rho_H,\sigma_H)$, set $\sigma:=\sigma_H/2$. Now consider the function $\Omega(\rho):=\sup_{I \in \ml{G}_{\rho}}|\omega_C(I)|$, clearly $\Omega(0)=0$. From now on we shall suppose that $\rho$ satisfies the following condition    
\beq{eq:smallnessrho}
4 r N \Omega(\rho) \leq a \mx{.} 
\eeq 
The analyticity\footnote{Obviously, $\Omega(\rho) \equiv 0$ for all $\rho$ in the case of an isochronous system, so that (\ref{eq:smallnessrho}) would impose no restrictions on $\rho$.} of $h(I)$ implies the existence of $C_h \in [1,+\infty)$ such that the value of $\rho$ can be determined as
\beq{eq:rhodet}
\rho:=\min\{\rho_H,a(4 r N C_h)^{-1}\} \mx{,}
\eeq
once $r$ and $N$ will be chosen.\\The scheme is constructed in such a way one can set $(\tilde{\rho}_*,\tilde{\sigma}_*):=(1-d)(\rho,\sigma)$.\\
As a consequence of Hypothesis \ref{hyp} and of the standard properties of analytic functions, one has
\beq{eq:decayhs}
\norm{H_m}{\rho,\sigma} \leq \ml{F} h^{m-1} e^{-a t},\qquad m \geq 1 \mx{,}
\eeq 
with $\ml{F}:=\ep \tilde{\ml{F}}$, where (see \cite[Lemma 5.2]{gio03}) $\tilde{\ml{F}}:=[(1+\exp(-\sigma/2))/(1-\exp(-\sigma/2))]^n$ and  
\beq{eq:h}
h:=\exp(-N\sigma/2)\mx{.}
\eeq 
\begin{lem}\label{lem:homologicaleq} Suppose that $\anorm{\Psi_s}{s} \leq M^{(s)} \exp(-a_s t)$, for some $M^{(s)}>0$. Then the solution of (\ref{eq:homological}) satisfies
\beq{eq:boundschi}
4a \anorm{\chi^{(s)}}{s+\ah} , 
4 \anorm{\pl_t \chi^{(s)}}{s+\ah}
\leq C_r M^{(s)} e^{-a_{s+1}t} \mx{,}
\eeq
where
$C_r:=2^{2n+4}(r/d)^n$.
\end{lem}
\proof Use (\ref{eq:homeqan}). Similarly to Lemma \ref{lem:homeq}, we choose $c_k^{(s)}(I,0):=-\int_{\RR^+} 
\exp(-(\omega(I)\cdot k)\tau) \psi_k^{(s)}(I,\tau)d \tau$. Note that $|c_k^{(s)}(I,0)| \leq M^{(s)} \exp(-(1-d_s)|k|\sigma) \int_{\RR^+} \exp(|\omega_C(I)||k|-a)\tau) d \tau < +\infty$ on $\Lambda_s$ by (\ref{eq:smallnessrho}). By using again (\ref{eq:smallnessrho}) one gets
\beq{eq:similarto}
|c_k^{(s)}(I,t)| \leq M^{(s)} e^{-(1-d_s)|k|\sigma} e^{\frac{a_s s}{4r}t} \int_t^{\infty} e^{a_s \left(\frac{s-4r}{4r} \right) \tau} d \tau \leq \frac{4}{a} M^{(s)} e^{-(1-d_s)|k|\sigma}  e^{-a_s \left(1-\frac{s}{2r}\right)t} \mx{.}
\eeq
The first of (\ref{eq:boundschi}) is easily recognised\footnote{Use the inequality $\sum_{|k| \geq (s-1)N} \exp(-\delta|k|\sigma) \leq \exp(-N \delta n (s-1) \sigma) (\sum_{m=0}^{+\infty}\exp(-\delta m \sigma))^n \leq (2/\delta)^n$, where in this case $\delta:=d_{s+\frac{1}{2}}-d_s= d/(2r)$. }
by (\ref{eq:seqa}). The second of (\ref{eq:boundschi}) follow from (\ref{eq:similarto}) and from (\ref{eq:homeqan}).  
\endproof
 
\begin{lem}\label{lem:numerical} Let $A,\Gamma,\tau > 0$ and consider the real-valued sequences $\{\kappa_s\}_{s \geq 1}$ and $\{\gamma_l\}_{l \geq 0}$ defined as 
\beq{eq:ausseq}
\kappa_l:=A \tau^{l-1} + \Gamma \sum_{j=1}^{s-1} \tau^{j-1} \kappa_{l-j},\qquad
\gamma_l := \Gamma \sum_{j=1}^l \tau^{j-1} \gamma_{l-j} \mx{,}
\eeq
where $\kappa_1$ and $\gamma_0$ are given. Define $\Delta:=\tau+\Gamma$, then for all $s \geq 2$ and $ l \geq 1$
\beq{eq:ricorrbound}
\kappa_s = (\Gamma \kappa_1+\tau A)\Delta^{s-2}, \qquad 
\gamma_l = \gamma_0 \Gamma \Delta^{l-1} \mx{.}
\eeq
\end{lem}
\proof We shall denote with (\ref{eq:ausseq}a) and (\ref{eq:ausseq}b) the first and the second of (\ref{eq:ausseq}), respectively. The same for (\ref{eq:ricorrbound}).
Let us suppose for a moment that (\ref{eq:ricorrbound}a) is proven, then choose $A=\Gamma \gamma_0$ and $\kappa_1=\Gamma \gamma_0 = \gamma_1$. By substituting in (\ref{eq:ricorrbound}a) one immediately gets (\ref{eq:ricorrbound}b). Hence we need only to prove (\ref{eq:ricorrbound}a).\\ For this purpose we use the well-known \emph{generating function} method (see e.g. \cite{wil}). Namely, define $g(z):=\sum_{n=1}^{\infty} w_n z^n$, multiply each equation obtained from (\ref{eq:ausseq}a) by $z^s$ as $s$ varies, then ``sum'' all the equations. This leads to $g(z)=[1-\Delta z]^{-1}(\kappa_1 (z-\tau z^2)+A \tau z^2)=(1+\Delta z + \Delta^2 z^2+\ldots)(\kappa_1 (z-\tau z^2)+A \tau z^2)=\kappa_1 z + (\Gamma \kappa_1 + \tau A) \sum_{n \geq 2} \Delta^{n-2} z^n$, which is the (\ref{eq:ricorrbound}a).
\endproof
%%%%%%%%%%%%%%%%%%%%%%%%%%%%%%%%%%%%%%%%%%%%%%%%%%%%%%%%%%%%%%
\subsection{Bounds on the generating function}
\begin{prop}\label{prop:psi} For all $s \leq r$, the following estimate holds
\beq{eq:indchi}
\anorm{\chi_s}{s+1/2} \leq (4a)^{-1} C_r \beta_s \ml{F} e^{-a_{s+1} t} \mx{,}
\eeq
where the sequence $\{\beta_s\}_{s=1,\ldots,r} \in \RR^+$ is determined by the following system  
\beq{eq:recsys}
\left\{
\begin{array}{rcl}
\beta_s & = & \de h^{s-1} + \frac{\Gamma}{s} \sum_{j=1}^{s-1} j \theta_{s-j} \\
\theta_l & = & \de \frac{\Gamma}{l} \sum_{j=1}^l j \beta_j \theta_{l-j} 
\end{array}
\right.
\eeq
with $\{\theta_l\}_{l=0,\ldots,r-1} \in \RR^+$ and  
\beq{eq:gamma}
\Gamma:=16 n r^2 C_r \ml{F} (a d^2 \rho \sigma )^{-1} \mx{,}
\eeq
under the conditions\footnote{From a ``computational'' point of view, first compute $\theta_1$ then proceed with $\beta_s,\theta_s$ for all $s=2,\ldots,r$.}
$\beta_1=\theta_0=1$.
\end{prop}
First of all note that by (\ref{eq:lietransform}) and (\ref{eq:decayhs}),  one has 
$\anorm{\Psi_1}{1} \leq \ml{F}\exp(-a_1 t)$ and $\anorm{E_0 H_m}{} \leq \ml{F} h^{m-1} \exp(-a_1 t)$ (recall (\ref{eq:seqa})). Hence,  given by $s \leq r$, we can suppose by induction to know $\beta_1,\ldots,\beta_{s-1}$ and $\tilde{\theta}_{0,m},\ldots,\tilde{\theta}_{s-2,m}$, for all $m \geq 1$, with $\beta_1=1$ and $\tilde{\theta}_{0,m}=h^{m-1}$, such that the the following bounds hold for all $j=1,\ldots,s-1$ and $l=0,\ldots,s-2$
\begin{subequations}
\begin{align}
\anorm{\Psi_j}{j} & \leq \beta_j \ml{F} e^{-a_j t} \mx{,}
\label{eq:indpsi} \\ 
\anorm{E_l H_m}{l+1} & \leq \tilde{\theta}_{l,m} \ml{F} e^{-a_{l+1}t} \label{eq:elhm} \mx{,} 
\end{align}
\end{subequations} 
By (\ref{eq:indpsi}) and Lemma \ref{lem:homologicaleq}, the bound (\ref{eq:indchi}) holds with $j$ in place of $s$. Hence by Prop. \ref{prop:poisson} with $G=\chi^{(j)}$, $F=E_{s-j-1}H_m$ then $\hat{d}=\max_{j=1,\ldots,s-1}\{d_{j+\ah},d_{s-j}\}=d_{s-\ah}$ and finally $\tilde{d}:=d_s-d_{s-\ah}=d/(2r)$, one has (by setting $\delta=0$)
\beq{eq:estone}
\begin{array}{rcl}
\anorm{\ml{L}_{\chi^{(j)}} E_{s-j-1} H_m}{s} & \leq & 8 r^2 (e d^2 \rho \sigma )^{-1}  \anorm{\chi^{(j)}}{j+1/2} \anorm{E_{l-j} H_0}{l-j+1/2}\\
& \leq & \Gamma \ml{F} \beta_j \gamma_{l-j} e^{-a_{l+1}t}
\end{array}
\eeq
where the property $a_{j+1} + a_{l-j+1} \geq a_{l+1}$ has been used. Recalling (\ref{eq:lietransform}), we have that (\ref{eq:elhm}) holds also for $l=s-1$, where
\beq{eq:thetal}
\tilde{\theta}_{l,m} = \frac{\Gamma}{l} \sum_{j=1}^l j \beta_j \tilde{\theta}_{l-j,m} \mx{.}
\eeq
Furthermore, it is easy to show from the latter that $\tilde{\theta}_{l,m}=h^{m-1} \tilde{\theta}_{l,1}$ in such a way, defined $\theta_l:=\tilde{\theta}_{l,1}$ one gets $\tilde{\theta}_{l,m}=h^{m-1}\theta_l$, and then the second of (\ref{eq:recsys}), provided $\theta_0=1$.
In conclusion, by using (\ref{eq:decayhs}), 
and the second of (\ref{eq:recsys}) in the definition of $\Psi_s$ as in (\ref{eq:homological}), we get that (\ref{eq:indpsi}) is satisfied if $\beta_s$ is defined as in the first of (\ref{eq:recsys}). Bound (\ref{eq:indchi}) follows from Lemma \ref{lem:homologicaleq}.
\endproof
\begin{prop}\label{prop:estbeta}
The sequence $\beta_s$ defined by (\ref{eq:recsys}) satisfies 
\beq{eq:betasric}
\beta_s \leq \tau^{s-1}/s \mx{,}
\eeq
for $s=1,\ldots,r$, if 
\beq{eq:choice}
\tau:=eh, \qquad \Gamma \leq h/(2r^2)  \mx{.}
\eeq
\end{prop}
\proof
The property (\ref{eq:betasric}) is trivially true for $s=1$, hence let us suppose it for $j=1,\ldots,s-1$ and proceed by induction with $\tau$ to be determined. Define $\tilde{\theta}_l:=\theta_l(\beta_j)|_{\beta_j=\tau^{j-1}/j}$, then $\hat{\theta}_l:=\tilde{\theta}_l/l$, obtaining $
\hat{\theta}_l=\Gamma \sum_{j=1}^{l} \tau^{j-1} \hat{\theta}_{l-j} 
$. Clearly $\theta_l \leq \tilde{\theta}_l \leq \hat{\theta}_l/l$, furthermore $\theta_0 = \tilde{\theta}_0=\hat{\theta}_0=1$. Hence, by Lemma \ref{lem:numerical} we have   
\beq{eq:boundtheta}
\theta_l \leq \Gamma \Delta^{l-1}/l \mx{.}
\eeq
Now choose $\tau,\Gamma$ as in (\ref{eq:choice}). By using (\ref{eq:decayhs}) and (\ref{eq:boundtheta}) in the first of (\ref{eq:recsys}) one gets that (\ref{eq:betasric}) is satisfied simply by checking that the inequality
\beq{eq:ineq}
y(s):=s+\frac{(s-1)}{2 r^2} \left(e+ \frac{1}{2 r^2} \right)^{s-1} \leq e^{s-1} 
\eeq
holds true for all\footnote{Clearly (\ref{eq:ineq}) holds for $s \leq r$ if $y(r)\leq \exp(r-1)$ for all $r \geq 3$ (let it be directly checked for $r=1,2$). Hence set $r=n+1$ and prove that $y(r)_{r=n+1} \leq \exp(n)$ for all $n \geq 2$, conclusion that is immediate as one can find that $y(n) \leq n+1 + 3 e^n/(4n)$.} $s=1,\ldots,r$.
\endproof
%%%%%%%%%%%%%%%%%%%%%%%%%%%%%%%%%%%%%%%%%%%%%%%%%%%%%%%%%%%%
\subsection{Estimates on the coordinates transformation}
From now on we shall suppose that $h$ and $\ep$ are chosen in such a way  
\begin{subequations}
\begin{align}
8 e h & \leq 1 \label{eq:finone}\\
2 r^ 2 \Gamma & \leq \sqrt{\ep} h  \label{eq:fintwo}
\end{align}
\end{subequations}
In particular, by definition and by (\ref{eq:choice}), this immediately implies that 
\beq{eq:smallnessone}
4 \Delta \leq 1
\eeq
As in \cite{gio03} it is used that, despite the generating sequence is finite, one can use the bound obtained from \ref{prop:psi} 
\beq{eq:chigen}
\anorm{\chi^{(s)}}{}\leq (4a)^{-1} C_r \ml{F} \beta_s e^{-a_{r+1}t} \mx{,} 
\eeq   
with $\beta_s$ satisfying (\ref{eq:betasric}) \tf{for all} $s$, as it would be, trivially, $\beta_{>r}=0$.
\begin{prop} Define $(I^{(r)},\ph^{(r)},\eta^{(r)}):=T_{\chi^{[r]}}(I,\ph,\eta)$. Then the following estimates hold
\beq{eq:boundstransform}
\anorm{I-I^{(r)}}{}, 
\anorm{\eta-\eta^{(r)}}{} \leq \frac{d \rho}{8}e^{-a_{r+1}t},\quad 
\anorm{\ph-\ph^{(r)}}{} \leq \frac{d \sigma}{8}e^{-a_{r+1}t} \mx{.}
\eeq
\end{prop}
\proof Let us start from the variable $I$. Firstly, note that $\anorm{I-T_{\chi^{[r]}}I}{}\leq \sum_{s \geq 1} \anorm{E_s I}{}$. In addition 
\[
\anorm{E_1 I}{2} = \anorm{\pl_{\ph}\chi^{(1)}}{2} \leq 2 n r (e d \sigma)^{-1} \anorm{\chi^{(1)}}{3/2} \leq D_{\sigma} \ml{F} \exp(-a_{r+1}t) \mx{,}
\] 
with $D_{\sigma}:=n r C_r/(2 d \sigma a)$ by Prop. \ref{prop:psi}. Hence suppose $\anorm{E_l I}{l+1} \leq \ml{F} u_l \exp(-a_{r+1}t)$ for all $l=1,\ldots,s-1$ with $u_1=D_{\sigma}$ and proceed by induction.\\
The bound of $E_l I$  can be treated in the same way of (\ref{eq:elhm}) with the difference that in this case the term $\lie{l} I$ appearing in $E_l I$ needs to be bounded separately by using (\ref{eq:indchi}) and a Cauchy estimate. This leads to $u_l=\beta_l D_{\sigma}+ \Gamma/l \sum_{j=1}^{l-1} j \beta_j u_{l-j}$. By using the same procedure used in the proof of Prop. \ref{prop:estbeta} for $\theta_l$ one gets $u_l \leq (D_\sigma/l) \Delta^{l-1}$. The required bound easily follows as $ \ml{F} \sum_{s \geq 1 } u_s \leq 2 \ml{F} D_{\sigma} \leq \Gamma d \rho \leq \sqrt{\ep} d \rho /8$, where the second inequality follows from (\ref{eq:smallnessone}) and the last one from (\ref{eq:fintwo}) then from (\ref{eq:finone}). The procedure for the variables $\ph$ and $\eta$ is similar. The analyticity of the transformation $\ml{N}_r:=T_{\chi^{[r]}}^{-1}$ easily follows from the bounds (\ref{eq:boundstransform}) and the invertibility of the Lie transform operator, see \cite{gio03}.
\endproof

\subsection{Bound on the remainder}
\begin{prop}\label{prop:remainder} Define $A:=10 \tilde{\ml{F}}$ then for all $r \geq 1$
\beq{eq:estrem}
\norm{\ml{R}^{(r+1)}}{(1-2d)(\rho,\sigma)} \leq \ep A e^{-(r+a_{r+1}t)} \mx{.} 
\eeq
\end{prop}
\proof
Define $(\rho',\sigma'):=(1-d)(\rho,\sigma)$. Now recall (\ref{eq:remaindernekho}) and suppose by induction, for all $l=1,\ldots,s-1$, $m=0,\ldots,s-2$ with $s \in \NN$ 
\beq{eq:ricorrem}
\norm{E_l H_0}{(1-(l/s)d)(\rho',\sigma')} \leq \ml{F} \epsilon_l \exp(-a_{r+1}t), \qquad  
\norm{E_m H_n}{(1-(m/s)d)(\rho',\sigma')} \leq \ml{F} \zeta_{m,n} \exp(-a_{r+1}t) \mx{.}  
\eeq
Indeed one can set $\zeta_{0,n}=h^{n-1}$ and $\epsilon_1=\beta_1=1$ as $\lie{1}H_0=-\Psi_1$ by (\ref{eq:homological}). We stress that, despite based on the same computations, the argument is conceptually different from the previous estimates as $s \in (r,+\infty)$ and the use of  $\delta$ in (\ref{eq:lie}) plays here a key role. More precisely, use Prop. \ref{prop:poisson} with $G=\chi^{(j)}$ and $F=E_{s-j}H_0$ hence $d''=0$ then $\hat{d}=d'=\delta=d(s-j)/s$ from which $\tilde{d}=(j/s)d$. This leads to $\norm{\lie{j} E_{s-j} H_0}{(1-d)(\rho',\sigma')} \leq \Gamma (s/j)\beta_j \epsilon_{s-j} \exp(-a_{r+1}t)$, implying\footnote{The use of (\ref{eq:lie}) with $\delta=0$ would have given $(s/j)^2$ instead of $(s/j)$, producing in this way a troublesome factorial in the estimates.} that the first of (\ref{eq:ricorrem}) holds for $l=s$ provided $\epsilon_s=\beta_s+\Gamma \sum_{j=1}^{s-1} \beta_j \epsilon_{s-j} = \Delta^{s-1}$, the latter by Lemma \ref{lem:numerical}. This implies $\norm{\sum_{l=1}^s E_{s-l}H_l}{(1-d)(\rho',\sigma')} \leq \ml{F} (s+1) \Delta^{s-1} \exp(-a_{r+1}t)$ by using (\ref{eq:decayhs}) and the trivial bound $h \leq \Delta$.
Similarly one finds $\zeta_{s,n}=h^{n-1}\Delta^{s-1}$, hence 
\[
(\ml{F} e^{-a_{r+1}t})^{-1} \ml{R}^{(r+1)} \leq \sum_{s > r} (2+s) \Delta^{s-1} = 
 \Delta^r \left( \frac{r+3}{1-\Delta}+ \frac{1}{1-\Delta^2} \right) \leq 2 (r+4) \Delta^r \mx{,}
\]
by (\ref{eq:smallnessone}). Noticing that $\ml{D}_{(1-2d)(\rho,\sigma)} \subset \ml{D}_{(1-d)^2(\rho,\sigma)}$, the bound (\ref{eq:estrem}) easily follows from (\ref{eq:smallnessone}) and from the simple inequality $(r+4)e^r \leq 5 (4^r)$.
\endproof
%%%%%%%%%%%%%%%%%%%%%%%%%%%%%%%%%%%%%%%%%%%%%%%%%%%%%%
\subsection{Parameters choice and perpetual stability} 
Let us discuss a possible choice of the parameters in such a way the convergence conditions are satisfied. More precisely 
by (\ref{eq:h}), condition (\ref{eq:finone}) holds if $N=\lceil 2 \sigma^{-1}(1+3 \log 2) \rceil$, where $\lceil \cdot \rceil$ denotes the rounding to the greater integer.  This implies that $ h \geq 1/(16 e) $, hence (\ref{eq:fintwo}) holds if $2^5 e r^2 \Gamma \leq \sqrt{\ep}$. Hence, recalling  (\ref{eq:smallnessrho}) and (\ref{eq:gamma}), this condition is achieved by choosing (see also \cite{gg85})
\beq{eq:r}
r:=\left\lfloor \left( \frac{\ep_a^*}{\ep} \right)^{\frac{1}{2 \gamma}}   \right \rfloor, \qquad \sqrt{\ep_a^*}:=\frac{a^2 d^{n+2} \rho_H \sigma^2}{2^{2n+19} e n C_h \tilde{\ml{F}} } \mx{,}
\eeq
where\footnote{Note that the threshold $\ep_a^*$ takes into account of the condition (\ref{eq:rhodet}) as we have used the obvious lower bound $\rho \geq a\rho_H(4 r N C_h)^{-1}$, immediate from (\ref{eq:rhodet}).} $\gamma=5+n$ and $\lfloor \cdot \rfloor$ denotes the rounding to the lower integer. The condition $\ep \leq \ep_a^*$, as in the statement of Theorem \ref{thmtwo}, clearly ensures that $ r \geq 1$. The final value of $\rho$ is determined with (\ref{eq:rhodet}).\\
Let us write the usual bound $|I(t)-I(0)| \leq |I(t)-I^{(r)}(t)|+|I^{(r)}(t)-I^{(r)} (0)|+|I^{(r)}(0)-I(0)|$. The first and third term of the r.h.s. are bounded by $\sqrt{\ep}d \rho/8$ by (\ref{eq:boundstransform}). As for the second one, from the equations of motion $\dot{I}^{(r)}=-\pl_{\ph}H^{(r)}=-\pl_{\ph}\ml{R}^{(r+1)}$, furthermore $\norm{\pl_{\ph} \ml{R}^{(r+1)}}{(1-2d)(\rho,\sigma)} \leq \ep A (e d \sigma)^{-1} \exp(-(r+a_{r+1}t))$ by a Cauchy estimate and by (\ref{eq:estrem}). Hence 
\beq{eq:lastineq}
|I^{(r)}(t)-I^{(r)} (0)| \leq \ep A (e d \sigma)^{-1} e^{-r} \int_0^t e^{-a_{r+1}s}ds \leq \ep A (a d e \sigma)^{-1} (2/e)^{r} \mx{,}
\eeq
as $a_{r+1}=a(2r-1)(2r-2)\ldots(r)/(2r)^r>a2^{-r}$. 
\begin{rem}
The bound (\ref{eq:lastineq}) is the key element beyond the perpetual stability, despite a normal form of finite order. The remainder, which is bounded by a constant in the classical Nekhoroshev estimate and then produces a linearly growing bound for the quantity $|I^{(r)}(t)-I^{(r)} (0)|$, is now summable over $\RR^+$. Hence, a restriction to exponentially large times is no longer necessary.
\end{rem}
It is immediate from (\ref{eq:lastineq}) that for all $\ep \leq \ep_a^*$ one has  $|I^{(r)}(t)-I^{(r)} (0)|\leq 2 \ep_a^* A (a d e^2 \sigma)^{-1}$ which is clearly smaller than $\sqrt{\ep} d \rho/4$ by (\ref{eq:r}). Hence $|I(t)-I(0)| \leq \sqrt{\ep} d \rho/2$.

\subsection*{Acknowledgements} The first author is grateful to Proff. D. Bambusi, L. Biasco, A. Giorgilli and T. Penati for very useful discussions on a preliminary version of this paper. 

\newcommand{\sort}[1]{}

\end{document}